\documentclass{elsart}

\usepackage{amsfonts}

\usepackage{amsmath}

\usepackage{amssymb}

\input{xy}
\xyoption{all}

\newcommand{\Abar}{\overline{A}}
\newcommand{\Ibar}{\overline{I}}
\newcommand{\Pbar}{\overline{P}}
\newcommand{\Ahat}{\hat{A}}
\newcommand{\QA}{\mathcal{Q}_A}
\newcommand{\QAbar}{\mathcal{Q}_{\overline{A}}}
\newcommand{\smod}{\underline{\textup{mod}}\,}
\newcommand{\rad}{\textup{rad}\,}
\newcommand{\Hom}{\textup{Hom}\,}
\newcommand{\sHom}{\underline{\textup{Hom}}\,}
\newcommand{\Ext}{\textup{Ext}^1}

\newcommand{\ind}{\textup{ind}\,}
\newcommand{\add}{\textup{add}\,}
\newcommand{\soc}{\textup{soc}\,}

\newcommand{\pd}{\textup{pd}\,}
\newcommand{\isomorphe}{\cong}
\newcommand{\CA}{\mathcal{C}_A}
\newcommand{\LA}{\mathcal{L}_{\overline{A}}}
\newcommand{\RA}{\mathcal{R}_{\overline{A}}}
\newcommand{\LC}{\mathcal{L}_{{C}}}
\newcommand{\RC}{\mathcal{R}_{{C}}}

\newcommand{\DA}{\mathcal{D}^b(\textup{mod}\, A)}
\newcommand{\bimodule}{_A\!DA _A}

\newcommand{\epito}{\twoheadrightarrow}
\newcommand{\monoto}{\hookrightarrow}

\newcommand{\comment}[1]{}

\newcommand{\zG}{\Gamma}
\newcommand{\zS}{\Sigma}
\newcommand{\zO}{\Omega}
\newcommand{\za}{\alpha}
\newcommand{\zb}{\beta}
\newcommand{\zg}{\gamma}
\newcommand{\zl}{\lambda}
\newcommand{\zs}{\sigma}

\begin{document}

\begin{frontmatter} 
  \title{Cluster categories and duplicated algebras} 
  \author{I. Assem\thanksref{nserc}}, \address{D\'epartement de Math\'ematiques,
  Universit\'e de Sherbrooke, Sherbrooke (Qu\'ebec), J1K 2R1, Canada}
  \ead{ibrahim.assem@usherbrooke.ca}
  \author{T. Br\"ustle\thanksref{nserc+}}, \address{D\'epartement de Math\'ematiques,
  Universit\'e de Sherbrooke, Sherbrooke (Qu\'ebec), J1K 2R1, Canada \emph{and} 
  Department of Mathematics, Bishop's University,  Lennoxville, (Qu\'ebec),
  J1M 1Z7, Canada} 
  \ead{thomas.brustle@usherbrooke.ca}
  \author{R. Schiffler\thanksref{nateq}}, \address{Department of Mathematics and
  Statistics, University of Massachusetts at Amherst, Amherst, MA
  01003-9305, USA}
  \ead{schiffler@math.umass.edu} 
  \author{G. Todorov}, \address{Department of Mathematics,
  Northeastern University, Boston, MA 02115, USA} 
  \ead{todorov.neu.edu}
  \thanks[nserc]{Partially supported by the NSERC of Canada}
  \thanks[nserc+]{Partially supported by the NSERC of Canada and the universities of Sherbrooke and Bishop's}
  \thanks[nateq]{Partially supported by the NATEQ of Qu\'ebec, Canada}
  
  \begin{abstract} 
Let $A$ be a hereditary algebra. We construct a fundamental domain for
the cluster category $\CA$ inside the category of modules over the 
duplicated algebra $\Abar$ of $A$. We then prove that there exists a
bijection between the tilting objects in $\CA$ and the tilting
$\Abar$-modules all of whose non projective-injective indecomposable
summands lie in the left part of the module category of $\Abar$.
  \end{abstract} 
\end{frontmatter} 

\setcounter{section}{-1}

\begin{section}{Introduction}

Cluster categories were introduced in \cite{BMRRT}, and
for type $A_n$ also in \cite{CCS}, 
as a means for a better understanding of
the cluster algebras of Fomin and Zelevinsky
\cite{FZ1,FZ2}. The indecomposable objects (without self-extensions) in the
cluster category correspond to the cluster variables in the cluster
algebra and the  tilting objects in the cluster category to the
clusters in the cluster algebra.
Our objective in this note is to give an interpretation of the cluster
category and its tilting objects in terms of modules over a finite dimensional algebra.
 Indeed, let $A$
be a hereditary algebra over an algebraically closed field then, by
Happel's theorem \cite{H}, the derived category of bounded
complexes over the category $\textup{mod}\, A$ of finitely generated right
$A$-modules is equivalent to the stable module category over the
repetitive algebra $\Ahat$ of $A$ (in the sense of Hughes and
Waschb\"usch \cite{HW}). The algebra $\Ahat$ is infinite dimensional
but, in order to study the cluster category, it suffices to look at a
finite dimensional quotient of $\Ahat$, namely the duplicated algebra
$\Abar$ of $A$ defined and studied in \cite{A,ANS}. The resulting
embedding  of $\textup{mod}\, \Abar$ into $\textup{mod} \,\Ahat$ induces a functor
$\bar{\pi}$ from $\textup{mod} \,\Abar $ to the cluster category $\CA$ of
$A$. We prove that the functor $\bar{\pi}$ induces a one-to-one
correspondence between the indecomposable objects in the cluster
category  and the non projective-injective $\Abar$-modules lying in
the left part $\LA$ of $\textup{mod}\,\Abar$, in the sense of Happel,
Reiten and Smal\o\   \cite{HRS} (we then say that $\LA$ is an 
{\em exact fundamental
domain} for the functor $\bar{\pi}$). This opens the way to our main
result.
\begin{thm} 
Let $A$ be a hereditary algebra. There exists a one-to-one
correspondence between the multiplicity-free tilting objects in the
cluster category $\CA$ of $A$ and the multiplicity-free tilting
$\Abar$-modules such that all non projective-injective indecomposable
summands of $T$ lie in $\LA$.
\end{thm}
This correspondence is given explicitly as follows. Since any
indecomposable projective-injective $\Abar$-module is necessarily a
summand of $T$, then $T=T_0\oplus \overline P$, where $\overline P$ is
a uniquely determined
projective-injective $\Abar$-module and $T_0$ has no projective-injective
summands. If all the indecomposable summands of $T_0$ lie in $\LA$,
then $\bar\pi(T_0)$ is a tilting object in $\CA$ and conversely, any
tilting object in $\CA$ is of this form.

Since duplicated algebras appear as a perfect context to view (cluster-)tilting
objects as actual tilting \emph{modules}, we investigate these algebras
further. In particular we show that the simply-laced Dynkin case
corresponds to representation-finite duplicated algebras, which, in
addition, are simply connected. In this case several techniques are
known for computing the tilting modules, allowing us to find the
clusters in the corresponding cluster algebra.

We now describe the contents of our paper. After a brief preliminary
section, devoted to fixing  the notation and recalling the main facts we
shall be using, the second section contains a detailed description of
the left part $\LA$. In the third section, we prove that $\LA$ is an
exact fundamental domain for the natural functor and we
prove our main result in section four. Our final section is devoted to
deduce related properties of the duplicated algebra.

\end{section}

\begin{section}{Preliminaries}
\begin{subsection}{Notation.} Throughout this paper, we let $A$ denote
  a hereditary algebra over an algebraically closed field $k$. We
  denote by $\textup{mod} \,A$ the category of finitely generated right
  $A$-modules and by $\ind A$ a full subcategory whose objects are
  representatives of the isomorphism classes of
  indecomposable objects in $\textup{mod} \,A$. The derived category of bounded
  complexes over $\textup{mod}\, A$ will be denoted by $\DA$. 
  For a vertex $x$ in the quiver $\QA$ of $A$, we write $e_x$ for the
  corresponding primitive idempotent and $S_x,P_x,I_x$, respectively, for
  the corresponding simple, indecomposable projective and
  indecomposable injective $A$-module. The functor $D=\Hom_k(-,k)$ is
  the standard duality between $\textup{mod}\, A$ and $\textup{mod}\, A^{\textup{op}}$,
  and $\tau_A=D\,Tr$, $ \tau^{-1}_A=Tr\,D$ are the Auslander-Reiten
  translations in $\textup{mod}\, A$. We refer to \cite{ARS} for further facts
  about $\textup{mod}\, A$,  and to \cite{R} for the tilting theory of $\textup{mod}\, A$.
\end{subsection}
\begin{subsection}{The cluster category $\CA$.} The cluster category
  $\CA$ of $A$ is defined as follows. Let $F$ denote the endofunctor
  of $\DA$ defined as the composition $\tau^{-1}[1]$, where $\tau $ is
  the Auslander-Reiten translation  in $\DA$ and $[1]$ is the shift
  functor. Then $\CA$ is the quotient category $\DA/F$. Its objects
  are the $F$-orbits of objects in $\DA$ and the morphisms are given
  by
\begin{equation}\Hom_{\CA}(\tilde X, \tilde Y) = \oplus_{i\in \Zset}
 \  \Hom_{\DA}(F^iX,Y)\end{equation}
where $X$ and $Y$ are objects in $\DA$ and $\tilde X, \tilde Y$ are
  their respective $F$-orbits. It is shown in \cite{K} that $\CA$ is a
  triangulated category. 
  Furthermore, the canonical functor $\DA \to \CA$ is a functor of triangulated categories.
  We refer to \cite{BMRRT} for facts about the
  cluster category.
\end{subsection}

\begin{subsection}{The duplicated algebra $ \overline{A} $} 
The duplicated algebra of a hereditary algebra $A$ is the
 matrix algebra  
 \begin{equation}\overline{A}=\left[
\begin{array}{cc}
\ A\ &\ 0\ \\ DA & A 
\end{array}\right]=
\left\{\left[\left.
\begin{array}{cc}
\ a\ &\ 0\ \\ q&b 
\end{array}\right] \  \right|\  {a,b\in A ,\atop q \in DA}\right\}
\end{equation}
with the ordinary matrix addition and the multiplication induced by
the bimodule structure of $DA$. Writing $1$ for the identity of $A$,
and setting
\begin{equation}e=\left[{\begin{array}{cc}
\ 1\ &\ 0\ \\0&0\end{array}}\right]
\quad \textup{and} \quad 
e'=\left[\begin{array}{cc}
\ 0\ &\ 0\ \\0&1\end{array}\right],
 \end{equation}
we see that $\Abar$ contains two copies of $A$ given respectively by 
$e \Abar e$ and by $e'\Abar e'$. 
In order to distinguish between these 
 we denote the first one by
$A$ and the second one by $A'$. Accordingly, $\QA'$ denotes the quiver
of $A'$, $x'$ the vertex of $\QA'$ corresponding to $x\in (\QA)_0$,
and $e'_x $ the corresponding idempotent. Let 
$\overline{S}_x,\overline{P}_x,\overline{I}_x$
denote respectively the simple, indecomposable projective and
indecomposable injective  module in
$\textup{mod}\, \Abar$ corresponding to $x\in(\QA\cup\QA')_0$.

The ordinary quiver $\QAbar$ of $\Abar$ is constructed as follows.
It contains $\QA$ and $\QA'$ as  full convex connected subquivers 
 and every vertex of $\QAbar$ lies in either $\QA$ or
$\QA'$. There is an arrow $x'\to y$ whenever $\rad(e'_x\Abar
e_y)/ \rad^2(e'_x\Abar e_y) \ne 0$.
Observe that $e'_x\Abar e_y=D(e_y A e_x)$ and therefore, if $e_y A e_x
\ne 0$ then there is a non-zero path in $\QAbar$ from $x'$ to
$y$. Also, since $e'_x \Abar \isomorphe D(\Abar e_x)$, each
$\overline{I}_x=\overline{P}_{x'}$
is projective-injective having
$S_x$ as a socle and $S_{x'}$ as a top. On the other hand, each $\overline{P}_x$ has
its support lying in $\QA$ and is therefore equal to the projective
$A$-module $P_x$. Dually, $\overline{I}_{x'}$ has its support lying
completely in $\QA'$ and equals the injective $A'$-module $I_{x'}$. 
For facts about the duplicated algebra, we refer to \cite{A,ANS}. 
\end{subsection}

\begin{subsection}{The repetitive algebra $\Ahat$}\label{section}
For our purposes, another description of $\Abar$ is needed. The
repetitive algebra  $\Ahat$ of the hereditary algebra  $A$,  is the infinite
matrix algebra
\begin{equation}\Ahat  \quad = \quad \left[
  \begin{array}{cccccccccc}
  \ddots &&&\ 0\ \\
&\ A_{m-1}\ \\
&Q_m&\ A_m\ \\
&&Q_{m+1}&\ A_{m+1}\ \\
&\ 0\ &&&\ddots
  \end{array}
 \right] \end{equation}
where matrices have only finitely many non-zero coefficients, $A_m=A$
and $Q_m= \bimodule$ for all $m\in \Zset$, all the remaining
coefficients are zero and multiplication is induced from the canonical
isomorphisms $A\otimes_A DA\isomorphe  \bimodule \isomorphe DA\otimes_A
A$ and the zero morphism $DA\otimes_A DA \to 0$, see \cite{HW}. Then
$\Abar$ is identified to the 
quotient algebra of $\Ahat$ defined by the surjection 
\begin{equation}\Ahat \quad \to \quad \left[
  \begin{array}{cc}
    A_0&0\\ Q_1&A_1
  \end{array}\right].
  \end{equation}
This identification yields an embedding functor
$\textup{mod}\, \Abar
\monoto \textup{mod}\, \Ahat$.
Similarly, the canonical surjection $\Abar \to e\Abar e=A$
yields an embedding functor $\textup{mod}\, A \monoto \textup{mod}\,
\Abar$.
Our first objective will be to look more closely at these embeddings.
\end{subsection}
\end{section}
\begin{section}{The left part $\LA$ of the duplicated algebra $\Abar$}
\begin{subsection} {Definitions and a preparatory lemma}
Let $C$ be any finite dimensional $k$-algebra, and $M,N$ be two
indecomposable $C$-modules. A {\em path} from $M$ to $N$ in $\ind C$ is
a sequence of non-zero morphisms
\begin{equation}
  \label{eq*} M=M_0\stackrel{f_1\,}{\to} M_1 \stackrel{f_2\,}{\to}
  \cdots \stackrel{f_t\,}{\to} M_t =N
\end{equation}
with all $M_i$ in $\ind C$. We denote such a path by $M\leadsto N$ and
say that $M$ is a {\em predecessor} of $N$ (or that $N$ is a {\em successor}
of $M$). 
When each $f_{i}$ in (\ref{eq*}) is irreducible, we say that (\ref{eq*}) is a {\em
  path of irreducible morphisms}. A path (\ref{eq*}) of irreducible morphisms
is {\em sectional} if $\tau_C\  M_{i+1}\ne M_{i-1}$ for all $i$ with
$1\le i\le t$. A {\em refinement} of (\ref{eq*}) is a path in $\ind C$: 
\begin{equation}M=M'_0\stackrel{f_1'\,}{\to} M_1' \stackrel{f_2'\,}{\to} \cdots \stackrel{f_s'\,}{\to} M'_s =N\end{equation}
with $s\ge t$ such that there exists an order-preserving injection 
$\zs: \{1,\ldots,t-1\}\to \{1,\ldots , s-1\}$ verifying $M_{i} = M'_{\zs(i)}$  for all $i$ with 
$1\leq i\le t$.

A full subcategory $\mathcal{C}$ of $\ind{C}$ is called {\em convex} in
$\ind C$ if, for any path (\ref{eq*}) from $M$ to $N$ in $\ind C$,
with $M,N$ lying in $\mathcal{C}$, all the $M_i$ lie in $\mathcal{C}$.

Useful examples of convex subcategories arise from the standard
embeddings $\textup{mod}\, A \monoto \textup{mod}\, \Abar $ and $\textup{mod}\, \Abar \monoto \textup{mod}\,
\Ahat$, as seen in \ref{section} above. 
We have the following lemma (see \cite[2.5]{A}, \cite[3.4,\,3.5]{T} or
\cite[4.1]{Y}), which will be used quite often when
 considering $A$-modules as $\Abar$-modules or $\Ahat$-modules. 

\begin{lem}
  \label{lemma1} 
  a) The embeddings $ \textup{mod}\, A \monoto
  \textup{mod}\, \Abar$ and $\textup{mod}\, \Abar \monoto \textup{mod}\,
  \Ahat$ are full, exact and preserve  
  indecomposable modules, almost split sequences and irreducible
  morphisms.\\
b) Under these embeddings, $\ind A$ is a full convex subcategory of
  $\ind \Abar$, closed under predecessors, and $\ind \Abar$ is a full
  convex subcategory of $\ind \Ahat$.
\end{lem}
\end{subsection}

\begin{subsection}{The left part}
Let again $C$ be a finite dimensional algebra.
Following Happel, Reiten and
Smal\o\   \cite{HRS}, we define the {\em left part} $\LC$ of
$\textup{mod}\,C$ to be the full subcategory of $\textup{mod}\,C$
consisting of all indecomposable $C$-modules 
such that if $L\leadsto M$, then the projective dimension $\pd L$ of
$L$ is at most one. The right part $\RC$ is defined dually.

Our objective now is to compute the left part of the module category of the
duplicated algebra $\Abar$ of a hereditary algebra $A$. We start by
observing that, by Lemma \ref{lemma1}, the 
complete slice of the Auslander-Reiten quiver $\zG(\textup{mod}\, A)$
 of $A$  consisting of the indecomposable injective $A$-modules embeds fully
  inside $\zG(\textup{mod}\, \Abar)$. The sources in this slice are the
  injectives $I_a$ with $a$ a sink in $\QA$. For each sink $a$ in
  $\QA$, the injective $A$-module $I_a$ is the radical of the projective-injective
  $\Abar$-module $\Ibar_a=\Pbar_{a'}$. 
   
 We recall that for any  algebra $C$ and any $L$ in $\textup{mod}\, C$, $\pd L\le 1$
 if and only if $\Hom_{C}(DC,\tau_C L)=0$ \ (see \cite[IX.1.7,
 p.319]{ARS} or \cite[p.79]{R}).
  \begin{lem}\label{lemma2} 
Let $M$ be an indecomposable $\Abar$-module. Then:
    \begin{enumerate}
  \item[a)]  If $M$ belongs to $\ind A$, then $M\in \LA$ and
  $\tau^{-1}_{\Abar}M \in \LA$. 
 \item[b)] If $M$ does not belong to $\ind A$, then there exist a
    sink  $a\in (\QA)_0$ and a path $\overline{P}_{a'} \leadsto M$.
\end{enumerate}  
  \end{lem}
  \begin{pf}
a) Any $A$-module $M$ admits a projective resolution in $\textup{mod}\, A$ of
    the form
\begin{equation}0\to P_1 \to P_0\to M \to 0\end{equation}
with $P_0$ and $P_1$ projective $A$-modules, hence projective
$\Abar$-modules. Thus the projective dimension of 
$M$ as an $A$-module and also as an $\Abar$-module is at most one. 
 This shows that $\ind A \subset \LA$, because $\ind A$ is closed
 under predecessors. 
 
 To see that $\tau_{\Abar}^{-1}M$ is in $\LA$, notice that, since $M$ is in $\ind A$, 
 $\Hom_{\Abar}(\Ibar_{x},M)= 0$ 
 for all injective 
 $\Abar$-modules $\Ibar_{x}$. 
 So $\pd_{\Abar}(\tau_{\Abar}^{-1}M)\le 1$ by the above
 remark. Furthermore, any non-projective predecessor $L$ of
 $\tau_{\Abar}^{-1}M$ lies in  
 $\ind A \cup\tau_{\Abar}^{-1}(\ind A)$, hence $\pd L\le 1$.

 b) Assume now that $M$ is not in $\ind A$. Then there exists $b\in (\QA)_0$ 
 such that $\Hom_{\Abar}(\overline {P}_{b'}, M) \ne 0$. If $b$ is a
 sink, we are done. If  not, consider the projective $A$-module
 $P_{b}$. Let $S_{a}$ be a simple submodule of $P_{b}$. Note that
 $S_{a}$ is projective since $A$ is hereditary. Therefore $S_{a}=
 P_{a}$ and $a$ is a sink. Then $\Hom_{A}(P_{a},P_{b})\ne0$ implies
 $\Hom_{A}(I_{a},I_{b})\ne0$, which induces a non-zero morphism  
 $\overline{P}_{a'}=\overline{I}_{a}\to
 \overline{I}_{b}=\overline{P}_{b'}$ of $\Abar$-modules. This yields
 the required path $\overline{P}_{a'} \leadsto M$. 
\qed
   \end{pf}
 \end{subsection}
\begin{subsection}{A characterization of the modules in $\LA$}   
  Before stating the next proposition, we recall that, by \cite[1.6]{AC},
  $\LA$ consists of all  
  $M\in \ind\Abar$, such that, if there exists a path from an
  indecomposable injective module to $M$, then this path can be
  refined to a path of irreducible morphisms, and any such refinement
  is sectional. 
\begin{prop}\label{proposition3}
   An indecomposable $\Abar$-module $M$ is in $\LA$ if and only if,
   whenever there exists a path $\overline{P}_{a'} \leadsto M$, with
   $a$ a sink in $(\QA)_0$, this path can be refined to a path of
   irreducible morphisms, and each such refinement is sectional. 
  \end{prop}
  \begin{pf}
Since the necessity follows directly from the above statement, we only
prove the sufficiency. Assume that $M$ satisfies the stated
condition. In order to prove that $M\in\LA$, it suffices to show that, 
if there exists a path $\Ibar_{x}\leadsto M$, with $\Ibar_{x}$
injective in $\textup{mod}\, \Abar$, then this path can be refined to a path of
irreducible morphisms, and any such refinement is sectional. Since
$\Ibar_{x}$ is not an $A$-module, it follows from Lemma
\ref{lemma2} b), that there exist a sink $a$ in $\QA$ and a path
$\Pbar_{a'}\leadsto\Ibar_{x}$, giving a path
$\Pbar_{a'}\leadsto\Ibar_{x}\leadsto M$. The conclusion follows at
once. 
\qed
\end{pf}
\end{subsection}
\begin{subsection}{Ext-injectives in $\LA$}
We now characterize the Ext-injectives in the additive full
subcategory $\add \LA$
of $\textup{mod}\, \Abar$ generated by the left part.
We recall from \cite{AS}  that, if ${\mathcal A}$ is an additive full subcategory
of $\textup{mod}\, \Abar$, 
closed under extensions, then an indecomposable  module $M$ in ${\mathcal
  A}$  is called
an Ext-$injective$ $in$ $ \mathcal A$ if $\Ext_{\Abar}(\quad
,M)\vert_{\mathcal A} = 0$.  
It is known that $M$ is Ext-injective in $\add \LA$ if and only if
$\tau_{\Abar}^{-1} M$ is not in 
$\add\LA$ (see \cite[3.4]{AS}). 
We denote by $\zS$ the set of all
indecomposable Ext-injectives in $\add\LA$. The following corollary
says that $\LA= \ind A \cup \zS$.

\begin{cor}\label{cor1} 
The following are equivalent for an $\Abar$-module $M$:
\\
a) $M$ is in $\zS$. 
\\
b) $M$ is in $\LA$ and $M$ is not in $\ind A$.
\\
c) $M$ is in $\LA$ and there exist a sink $a\in(\QA)_0$ and a path
    $\overline{P}_{a'}\leadsto M$.
\\
d) There exist a sink $a\in(\QA)_0$ and a path
    $\overline{P}_{a'}\leadsto M$ and any such path is
refinable to a sectional path. 
\end{cor}
\begin{pf} 
\\ a) implies b) since $\ind A \cup \tau^{-1}(\ind A) \subset \LA$ by
Lemma \ref {lemma2} a). 
\\ b) implies c) follows from  Lemma \ref {lemma2} b).
\\ c) implies d) follows from Proposition \ref {proposition3}.
\\ d) implies a) Proposition \ref {proposition3} implies that $M$ is
in $\LA$. The fact that there exist a  sink $a\in(\QA)_0$ and a path 
    $\overline{I}_{a}=\overline{P}_{a'}\leadsto M$ 
    (hence a sectional path), implies that 
    $\Hom_{\Abar}(\overline{I}_{a},M) \ne 0$ 
    by \cite[III.2.4, p.239]{ARS}. By the remark before Lemma
    \ref{lemma2}, it follows that 
    $\pd_{\overline{A}}(\tau_{\Abar}^{-1}M) \ge 2$ and therefore
    $\tau_{\Abar}^{-1}M$ is not in $\LA$.
\qed
\end{pf}
\begin{cor}\label{cor2} 
  The set $\zS$ of all indecomposable Ext-injectives in
  $\add\LA$ consists of all the projective-injectives lying in $\LA$
  as well as all the modules of the form
  $\tau^{-1} I_x$ with $ x \in (\QA)_0$, that is
  \begin{equation}\zS = \{ \tau^{-1}_{\Abar} I_{x}\mid {x\in (Q_{A})_{0}}\} \cup
    \{\Pbar _{x'}\mid \Pbar_{x'}\in  \LA\}.
  \end{equation}
\end{cor}
\begin{pf}
Clearly, projective-injective modules which lie in $\LA$ belong to
$\zS$. Now let $x\in (\QA)_0$ and  consider  
$ \tau^{-1}_{\Abar} I_{x}$. Let $a$ be a sink and $I_{a}$ be a maximal
indecomposable injective $A$-module such that there is an epimorphism
$I_{a} \to I_{x}$. Then there is a non-zero map 
$\tau_{\Abar}^{-1}I_{a} \to \tau_{\Abar}^{-1}I_{x}$ and therefore a path 
$\Ibar_{a} \to \tau_{\Abar}^{-1}I_{a} \to
\tau_{\Abar}^{-1}I_{x}$. Since $\tau_{\Abar}^{-1}I_{x}$ is in $\LA$,
it follows that $\tau_{\Abar}^{-1}I_{x}$ is in $\zS$ by Corollary \ref{cor1} d).

Conversely, suppose $X$ belongs to $\zS$, but is not a
projective-injective lying in $\LA$. By Corollary \ref{cor1}, 
there exists a sink $a$ and a sectional path in $\ind \Abar$
\begin{equation}
  \Pbar_{a'}=\Ibar_{a}= M_{0}\to M_{1} \to \ldots \to M_{t} = X 
\end{equation}
with $t\ge 1$ and $M_1= \Ibar_a/ S_a = \tau^{-1}_{\Abar}\,I_a$.
We claim that no $M_i$ (with $i\ge 1$) is a projective
$\Abar$-module. Indeed, assume first that $M_i$ (with $i\ge 1$) is
projective-injective. By hypothesis, $i<t$. Then $M_{i-1}=\rad M_i$
and $M_{i+1}= M_i/\soc M_i = \tau^{-1}_{\Abar}\, M_{i-1}$,
contradicting the sectionality of the above path.
On the other hand, for any $i\le t$, $\Hom_{\Abar}(\Pbar_{a'},M_i)\ne
0$ hence $M_i$ is not an $A$-module, and a fortiori not projective in
$\textup{mod}\, A$. This establishes our claim. We infer the existence of a
sectional path in $\ind \Abar$
\begin{equation}
  I_a=\tau_{\Abar}\, M_{1} \to \tau_{\Abar}\, M_2 \ldots \to
  \tau_{\Abar}\, M_{t} = \tau_{\Abar}\, X .
\end{equation}
Since $X\in\LA$, then, for any $i\le t$, $M_i\in \LA$ and so $\pd
M_i\le 1$ implying that $\Hom_{\Abar}(\Pbar_{x'},\tau_{\Abar}\,
M_i)=0$ for any $x\in(\QA)_0$. This shows that the above  path lies
entirely in $\textup{mod}\, A$. Since $I_a$ is injective, all the modules on it
are injective. In particular, there exists $x\in(\QA)_0$ such that $
\tau_{\Abar}\, X= I_x$.
\qed
\end{pf}
We now give another expression for the set of all indecomposable
Ext-injectives in $\add\LA$. For this, we need to recall that, if $M$
is an $\Abar$-module, then its first cosyzygy $\zO_{\Abar}^{-1}\,M$ is the
cokernel of an injective envelope $M\to \Ibar$ in $\textup{mod}\,\Abar$.

\begin{prop}\label{prop1} 
Let $x\in (Q_{A})_{0}$. Then $\zO^{-1}_{\Abar} \,P_{x} \isomorphe
\tau^{-1}_{\Abar}  I_{x}$. 
Consequently, 
\begin{equation}
\zS= \{ \zO^{-1}_{\Abar} P_{x}\mid {x\in (Q_{A})_{0}}\} \cup
    \{\Pbar _{x'}\mid \Pbar_{x'}\in  \LA\}.
\end{equation}
\end{prop}
\begin{pf}
We prove this by induction on the Loewy length of the projective
module $P_{x}$. Recall that the Loewy length of a module $M$ is the
smallest integer $i$ with $\rad^i M =0$.
Let $P_{a}= S_{a}$ be a simple projective module. 
Then $\zO^{-1}_{\Abar} P_{a} \isomorphe \Ibar_{a}/ S_{a}$. On the
other hand, from  the almost split sequence: 
\begin{equation}
0 \to  I_{a} \to \Ibar_{a} \oplus I_{a}/S_{a} \to \Ibar_{a}/S_{a} \to
0
\end{equation}
it follows that $\zO_{\Abar}^{-1}\,P_{a} \isomorphe \tau_{\Abar}^{-1} I_{a}$ for any
sink $a$, which proves our claim in this case. 
For an indecomposable non-simple projective $P_{x}$ let the radical be
$\rad P_{x}= \oplus\, P_{y_{i}}$. Then there are the following isomorphisms
of the injective envelopes: 
$I_{0}\, (P_{x}) = I_{0}\,(\rad P_{x}) \isomorphe \oplus\,
I_{0}\,(P_{y_{i}})$. 
Then 
$\zO_{\Abar}^{-1}(P_{x}) = I_{0}\,(P_{x})/ P_{x}$ and 
$\zO_{\Abar}^{-1} (\rad P_{x}) = \oplus\, \zO_{\Abar}^{-1}(P_{y_{i}})\isomorphe
I_{0}\,(P_{x})\, /\,(\oplus \,P_{y_{i}})$. 
A simple application of the snake lemma yields
 $\zO_{\Abar}^{-1} (P_{x}) \isomorphe \zO_{\Abar}^{-1} (\rad P_{x})
\,/\,S_{x}$.
Now, it is easy to see that there is an almost split sequence
\begin{equation}
0 \to  I_{x} \to \left( \oplus\,\left(\tau^{-1}_{\Abar}
I_{y_{i}}\right)\right)\oplus 
I_{x}/S_{x} \to \tau^{-1}_{\Abar} I_{x} \to 0. 
\end{equation}
Since each morphism in this sequence is irreducible, it is  either a
monomorphism or  an epimorphism. 
Since $S_x$ is the kernel of the morphism $I_{x}\to I_{x}/S_{x}$,
another application of the snake lemma  and the induction hypothesis 
$\tau^{-1}_{\Abar} (I_{y_{i}})\ \isomorphe \ \zO_{\Abar}^{-1}(P_{y_{i}})$
yield 
\begin{equation} 
 \tau^{-1}_{\Abar} I_{x}\ \isomorphe \ \left.\left(
 \oplus\,\left(\tau^{-1}_{\Abar} 
 I_{y_{i}}\right)\right)\right/ S_{x} \ \isomorphe\ \left. \left(
 \oplus\,\left(\zO^{-1}_{\Abar} P_{y_{i}}\right)\right)\right/ 
 S_{x}\  \isomorphe\ \zO_{\Abar}^{-1}(P_{x}).
\end{equation}
 \qed
\end{pf}
\end{subsection}
\end{section}
\begin{section}{Fundamental domain for the cluster category.}
\begin{subsection}{$\LA$ as a subcategory of $\textup{mod}\, \Ahat$ }
As a consequence of the above description, the left part $\LA$ is nicely
embedded in $\textup{mod}\, \Abar$, and thus in $\textup{mod}\, \Ahat$. 
\begin{cor}
  The embedding $\LA \monoto \textup{mod}\, \Abar \monoto
  \textup{mod}\, \Ahat$ is full, 
  exact and preserves indecomposable modules, irreducible morphisms
  and almost split sequences.
\end{cor}
\end{subsection}

\begin{subsection}{Relation between $\LA$ and $\CA$}
We are now able to describe an exact fundamental domain for the
cluster category $\CA$ inside $\textup{mod}\,\Abar$, and  actually inside $\LA$. Indeed, since $A$
is hereditary, and thus of finite global dimension, we have a
triangulated equivalence $\DA\isomorphe \underline{\textup{mod}\,}\Ahat$
(see \cite{H}). Let
 \begin{equation}\hat \pi : \textup{mod}\,\, \Ahat \epito \underline{\textup{mod}\,}
\Ahat \isomorphe \DA \epito \CA\end{equation}
  be the canonical functor. We define an {\em exact fundamental domain} for
  $\hat \pi$ to be a full convex subcategory of $\ind
  \Ahat$ which contains exactly one point of each fibre $\hat \pi
  ^{-1}(X)$, with $X$ an indecomposable object in $ \CA$.

We recall at this point that $\ind \Abar$ is a full convex subcategory
of $\ind\Ahat$.

\begin{thm}\label{cor5} 
The functor $\hat\pi$ induces a one-to-one correspondence between
the non projective-injective modules in $\LA$ and the indecomposable
objects in $\CA$. In particular, $\LA$ is an exact fundamental domain for $\hat\pi$.
\end{thm}
\begin{pf} Since
  $\LA$ is a full  convex subcategory of $\ind \Abar$, it is also
  convex  inside $\ind \Ahat$. 
 Furthermore, the non projective-injective modules in $\LA$ are just
  the modules in  $\ind A$ and those of  $\{\zO_{\Abar}^{-1} P_x\mid x \in
  (\QA)_0\}$. The statement follows at once from the definition of
  $\CA$ and from the fact that under the triangle equivalence
  $\DA\isomorphe \underline{\textup{mod}}\, \Ahat$, the shift of $\DA$
  corresponds to $\zO_{\Ahat}^{-1}$ (see \cite{H}).
\qed
\end{pf}
\end{subsection}
\end{section}
\begin{section}{Tilting modules vs tilting objects}
\begin{subsection}{The main theorem} 
In this section, we prove our main theorem, which
  compares the tilting $\Abar$-modules with the tilting objects in
  $\CA$. For this purpose, we assume without loss of generality that
  our tilting modules and our tilting objects are
  multiplicity-free. We start by observing that, if $T$ is a tilting
  $\Abar$-module, then every indecomposable projective-injective
  $\Abar$-module is a  direct summand of $T$. Hence $T$ decomposes
  uniquely as $T=T_0\oplus e'\Abar$, where $T_0$ has no
  projective-injective direct summands. We say that $T$ is an
 $\mathcal{L}$-{\em tilting module} if $T_0\in\add\LA$.
 
 We denote by $\bar\pi: \textup{mod}\,\Abar \to \CA$, the composition of
 the inclusion  
 $\textup{mod}\, \Abar \monoto \textup{mod}\,\Ahat$ and the functor
 $\hat\pi$. By abuse of notation, the modules will be often denoted by
 the same letter even when considered as objects in different
 categories. 
   \begin{thm}\label{mainthm}
    There is a one-to-one correspondence
\begin{eqnarray}
      \left\{ \mathcal{L}-{tilting \ modules}\right\}
      &\longleftrightarrow& \left\{ {Tilting\  objects\  in
\       }\CA\right\}\nonumber\\
{given\  by } \hspace{2cm}      T=T_0\oplus e'\, \Abar&\longleftrightarrow& \bar\pi(T_0).\nonumber
\end{eqnarray}
 \end{thm}
  \begin{pf}
    Let  $T=T_0\oplus e'\Abar$ be an $\mathcal{L}$-tilting module and let
    $X=\bar\pi(T_0)$. Say $T_0=\oplus_{i=1}^n T_i$
    where the  $T_i$ are pairwise non-isomorphic indecomposable
    $\Abar$-modules.  Then $X=\oplus_{i=1}^n     X_i$ 
    with
    $X_i=\bar\pi(T_i)$.
    We first notice that, clearly, the number $n$ of indecomposable
    summands of $T_0$ is equal to the rank of the Grothendieck group
    of $A$. Hence, in order to show that $X$ is a tilting object in
    $\CA$, it suffices to prove that $\Ext_{\CA}(X,X)=0$.
   Suppose to the contrary that  there exist $i,j$ such
    that $\Ext_{\CA}(X_j,X_i)\ne 0$. Since $\Ext$ is symmetric in the
    cluster category by \cite[1.7]{BMRRT}, we also have
    $\Ext_{\CA}(X_i,X_j)\ne 0$. Thus 
    there are non-zero morphisms $X_i\to \tau_{\CA} X_j$ and  $X_j\to
    \tau_{\CA} X_i$ in $\CA$. Let
    $\hat F=\zO_{\Ahat}^{-1}\tau_{\Ahat}^{-1}$. Then there exist integers
    $s,t\ge 0$  such that the previous morphisms lift to non-zero morphisms in $\smod \Ahat$ 
\begin{equation}T_i\to \hat F^s\tau_{\Ahat} T_j \quad \textup{and} \quad  T_j\to  
\hat F^t\tau_{\Ahat} T_i,\end{equation}
by definition of the cluster category and the triangulated structure
of $\textup{mod}\,\Ahat$, see \cite{H}.  
    Moreover $s\ne 0$ and $t\ne 0$ since by hypothesis
    $\Ext_{\Abar}(T_j,T_i)=\Ext_{\Abar}(T_i,T_j)=0$. Now $T_i,T_j$ are
    in $\LA=\ind A \cup \zS$. We then have $3$ cases to consider.
    \begin{enumerate}
      \item $T_i,T_j\in \zS$. Then $X_i$ and $X_j$ lie on a slice
      of $\CA$, 
      hence $\Ext_{\CA}(X_i,X_j)=0$, a  contradiction.
      \item $T_i,T_j\in \ind A$. If $s=1$, then there is a non-zero
      morphism $T_i\to \hat F\tau_{\Ahat} T_j
      =\zO_{\Ahat}^{-1} T_j$ in  $\smod \Ahat$. But this is impossible since 
      \begin{equation} \sHom_{\Ahat}(T_i,\zO_{\Ahat}^{-1} T_j) =
      \Hom_{\DA}(T_i,T_j[1]) = 
      \Ext_A(T_i,T_j)=0\end{equation}  
      where we have identified the modules $T_i$ and $T_j$ with the
      corresponding stalk complexes in $\DA$.
      Assume thus that $s\ge 2$.
      Now, either $\tau_{\Ahat} T_j$ is an $A$-module, or $T_j$ is a
      projective $A$-module, and then  $\tau_{\Ahat} T_j[1]$ is an
      $A$-module. But this fact and the structure of the morphisms in
      the derived category (see \cite{H}) imply that 
      \begin{equation}\underline{\Hom}_{\Ahat}(T_i,\hat F^s\tau_{\Ahat} T_j) =
      \Hom_{\DA}(T_i, F^s\tau_{\Ahat} T_j) =0,\end{equation} 
      again a contradiction.
      \item $T_i\in\ind A, T_j\in\zS$. 
      Then by Proposition \ref{prop1}, there exists an indecomposable
      projective $A$-module $P_x$ such that $T_j=\zO^{-1}_{\Ahat}\,
      P_x$. Since $\Ahat$ is self-injective, it follows from
      \cite[IV.3.7]{ARS} that
      $\zO^{-2}_{\Ahat}=\nu_{\Ahat}\tau^{-1}_{\Ahat}$ where
      $\nu_{\Ahat}$ is the Nakayama functor in
      $\textup{mod}\,\Ahat$. Thus $\hat 
      F \,\tau_{\Ahat}\, T_j= \zO^{-2}_{\Ahat}\, P_x = \nu_{\Ahat} \,
      \tau^{-1}_{\Ahat}\, P_x$ which is an $A'$-module (unless $A$ is
      of Dynkin type $A_n$, linearly oriented
      and $P_x$ is projective-injective, in
      which case  $\hat F \,\tau_{\Ahat}\, T_j[-1] $ is an
      $A'$-module). Therefore the modules $T_i$ and  $\hat
      F^s \,\tau_{\Ahat}\, T_j$ have disjoint supports for any $s\ge
      1$. Therefore  
      $ \underline{\Hom}_{\Ahat}(T_i,\hat F^s\tau_{\Ahat} T_j) =0$ for
      any $s\ge 1$, contradiction.
    \end{enumerate}
  This completes the proof that $X=\hat\pi(T_0) $ is a tilting object
  in $\CA$.

Conversely, let $X=\oplus_{i=1}^n X_i$ be any tilting object in
$\CA$, where we assume that the objects $X_i$ are indecomposable and
pairwise non-isomorphic.
By Theorem \ref{cor5}, there exists, for each $i$ with $1\le i\le
n$,  a unique module $T_i\in\LA$ in the fibre $\hat\pi^{-1}(X_i)$.
Let  $T_0=\oplus_{i=1}^n T_i$. 
Then, clearly $\hat\pi(T_0)=X$.  
We want to show that $T=T_0\oplus e'\Abar$ is an $\mathcal{L}$-tilting
$\Abar$-module. Since $T_0\in\add\LA$ by construction
and, on the other hand, the number of indecomposable summands of $T_0$
is equal to the rank of the Grothendieck group of $A$, we only have to
prove that $\Ext_{\Abar}(T,T)=0$.
Suppose to the contrary, that there exist $i,j$ such that
$\Ext_{\Abar}(T_i,T_j)\ne 0$. Then $\Hom_{\Abar}(T_j,\tau_{\Abar}
T_i)\ne 0$.
In particular, $T_i $ is not projective in mod$\Abar$.
 Now, $T_i\in\LA$ implies
that $\tau_{\Abar} T_i=\tau_{\Ahat} T_i$. By Lemma \ref{lemma2} and
Corollary \ref{cor1}, we also have
$\tau_{\Abar} T_i\in\textup{ind}\, A$. Therefore $\Hom_{\Abar}(T_j,\tau_{\Abar}
T_i) \ne 0$ implies that $T_j\in\textup{ind}\, A$ (because $\ind A$ is
closed under predecessors in $\ind\Abar$).
Thus $\Hom_A(T_j,\tau_{\Abar} T_i)\ne 0$ and then
$\Ext_{\CA}(X_i,X_j)\ne 0$, contradiction. \qed
  \end{pf}

\end{subsection}
\begin{subsection}{Example}
Let $A$ be given by the quiver
$$
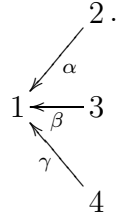
\xymatrix@C=20pt {&\ar[ld]^\za2\\
1&\ar[l]^\zb3\\
&\ar[lu]^\zg4}.
$$
 Then the ordinary quiver of $\Abar$ is given by  
$$
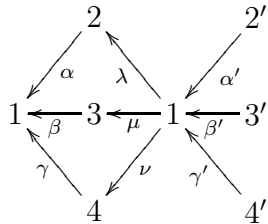
\xymatrix@C=20pt {&\ar[ld]^\za2&&\ar[ld]^{\za'} 2'\\
1&\ar[l]^\zb 3 &\ar[lu]^\zl\ar[l]^{\mu}\ar[ld]^{\nu} 1&\ar[l]^{\zb'}3'\\
&\ar[lu]^\zg4 &&\ar[lu]^{\zg'}4'}
$$
bound by the relations $\zl\za=\mu\zb=\nu\zg, \
\za'\mu=\za'\nu=\zb'\nu=\zb'\zl=\zg'\zl =\zg'\mu=0$. The
Auslander-Reiten quiver of $\Abar$ is given by

$$\xymatrix@C=10pt
{\save[0,0].[5,7] *[F-:]\frm{}
\\
&&&&&\circ\ar[rdd]&&&&&\circ\ar[rd]&&&&&&&\\
&\diamond\ar[rd]&& \bullet\ar[rd]&& \bullet\ar[rd]&& \bullet\ar[rd]\ &&
  \bullet\ar[rd]\ar[ru]&\circ\ar[rd]& \bullet\ar[rd]&& \bullet\ar[rd]&
  &\bullet \\
\ \bullet\ar[ru]\ar[r]\ar[rd]& \diamond\ar[r]&
\bullet\ar[ru]\ar[r]\ar[rd]& \bullet\ar[r]&
\bullet\ar[ru]\ar[r]\ar[rd]\ar[ruu]& \bullet\ar[r]&
\bullet\ar[ru]\ar[r]\ar[rd]& \bullet\ \ar[r]&
\bullet\ar[ru]\ar[r]\ar[rd]& \bullet\ar[r]\ar[ru]&
\bullet\ar[ru]\ar[r]\ar[rd]& \bullet\ar[r]&
\bullet\ar[ru]\ar[r]\ar[rd]& \bullet\ar[r]&
\bullet\ar[ru]\ar[r]\ar[rd]& \bullet\\ 
&\bullet\ar[ru]&& \diamond\ar[ru]&& \bullet\ar[ru]&&
    \diamond\ \ar[ru]&& \bullet\ar[ru]\ar[rd]&& \bullet\ar[ru]&&
	  \bullet\ar[ru]&& \bullet \\
&&&&&&&&&&\circ\ar[ru]&\\
&&&\LA\\
\restore
}
$$

where we have indicated the left part $\LA$. We have also indicated an
$\mathcal{L}$-tilting module $T=T_0\oplus e'\Abar$, where $T_0 \in
\add \LA$. The summands of $T_0$ are indicated by diamonds and the
(projective-injective) summands of $e'\Abar$ by circles.
\end{subsection}
\end{section}
\begin{section}{More on duplicated algebras of hereditary algebras}
It follows from our main theorem that the duplicated algebras of
hereditary algebras are quite a natural class to consider, since all
the tilting objects of the cluster category correspond to the actual
modules over the duplicated algebras. In this section we study other
properties of these algebras, which  are consequences of the description
of the left part $\LA$ and the $\textup{Ext}$-injectives as done in
the previous sections. 

We recall that a finite dimensional algebra $C$ is called {\em left (or
  right) supported} provided the class $\add\LC$ (or $\add\RC$) is
  contravariantly finite (or covariantly finite, respectively) in
  $\textup{mod}\,C$, see \cite{ACT,ACPT}. 
\begin{cor} The duplicated algebra 
  $\Abar$ of a hereditary algebra  is both left and right supported.
\end{cor}
\begin{pf}
  By \cite[3.3]{ACT}, the canonical module $T=U\oplus V$ (with
  $U=\oplus_{X\in\zS} X$ and $V=\oplus_{\Pbar_x\notin \LA} \,\Pbar_x$) is
  a partial tilting module. Now  the
  number of its indecomposable summands equals the number of
  isomorphism classes of indecomposable injective $A$-modules plus the
  number of 
  isomorphism classes of indecomposable projective-injective
  $\Abar$-modules. Hence $T$ is a tilting module and  $\Abar$
  is left supported, by
  \cite[thm. A]{ACT}. The other statement
  follows by symmetry.
\qed
\end{pf}
\begin{rem}
The assumption that $A$ is a hereditary algebra is essential. If $A$
is a tilted algebra which is the endomorphism algebra of a regular
tilting module, then it is easily seen that $\Abar$ is neither left
nor right supported.
\end{rem}

Equivalent statements to duplicated algebras being 
 representation-finite are given in the next corollary.
We recall that an algebra $C$ is said to be a {\em laura algebra}
\cite{AC} provided the class $\ind C \setminus (\LC\cup \RC)$ contains
only finitely many indecomposables.

\begin{cor}
  Let $A$ be a hereditary algebra. 
The following conditions are equivalent:
$$\begin{array}{ll}
(a) &  \Abar \ is\  a\  laura\  algebra.\\
(b) & A \  is \  of \ Dynkin\  type. \\
(c) &\Abar \  is \ \textup{{\it representation-finite}}.\\
\end{array}$$
If  this is the case, then $\Abar$ is simply connected.
\end{cor}
\begin{pf}
  We  denote by $\zS'$ the set of all indecomposable Ext-projectives in
  $\add\RA$.
By Lemma \ref{lemma2}, Corollary \ref{cor1} and their duals, the duplicated
  algebra $\Abar$ is laura if and
    only if the class
$[\tau_{\Abar}\,\zS, \tau^{-1}_{\Abar}\, \zS'] $
of all the 
$M\in \ind \Abar$ such that there exists a path $L\leadsto M\leadsto N,$
  with $\tau^{-1}_{\Abar}  L\in \zS$ and $\tau_{\Abar} N \in \zS'$
consists of finitely many indecomposables. Now, by \cite[2.6]{A} this
  class is an exact fundamental domain for the module category over
  the trivial extension $T(A)$ of $A$ by its minimal injective
  cogenerator $DA$. Therefore $\Abar$ is laura if and
    only if $T(A)$ is representation finite, or, by \cite{T}, if and
    only if $A$ is of Dynkin type which, by \cite[2.6]{A} is the case
  if and only if $\Abar$ is representation-finite. 
The last statement follows from \cite[2.7]{A}.
\qed
\end{pf}

\begin{rem}
  Assume $A$ to be representation-infinite. Then, of course, Theorems
 \ref{cor1} and \ref{mainthm} still apply. In this case as well, a
 good description of the module category of the duplicated algebra
 $\Abar$ is known (see \cite{A,ANS}) and, at least in the tame
 case, it is possible to compute explicitly the $\mathcal{L}$-tilting
 modules. 
\end{rem}

\begin{ack}
The authors wish to thank Raymundo Bautista for very useful
discussions. 
\end{ack}
\end{section}

\end{document}